\documentclass[11pt,a4paper]{amsart}

\usepackage[all]{xy}
\usepackage{amscd}

\newtheorem{theorem}{Theorem}[section]
\newtheorem{corol}[theorem]{Corollary}
\newtheorem{prop}[theorem]{Proposition}
\newtheorem{conj}[theorem]{Conjecture}
\theoremstyle{definition}
\newtheorem{defin}[theorem]{Definition}
\newtheorem{exam}[theorem]{Example}
\theoremstyle{remark}
\newtheorem*{erem}{Remark}

\def\iff{if and only if }
\def\oc{one-to-one correspondence }

\def\set#1{\left\{\,#1\,\right\}}
\def\setsuch#1#2{\left\{\,#1\,|\,#2\,\right\}}
\def\lst#1#2{ #1_1 , #1_2 , \dots , #1_{#2} }
\def\row#1#2{\left( #1_1 , #1_2 , \dots , #1_{#2} \right)}

\def\rad{\mathop\mathrm{rad}}
\def\Aut{\mathop\mathrm{Aut}\nolimits}
\def\ob{\mathop\mathrm{ob}\nolimits}
\def\Gr{\mathop\mathrm{Gr}\nolimits}
\def\supp{\mathop\mathrm{supp}\nolimits}
\def\pam{\mathop\mathrm{par}\nolimits}
\def\bp{{\bullet}}	\def\vec{\mathbf{vec}}

\def\md{\mbox{-}\mathbf{mod}}
\def\gdim{\mathop\mathrm{gl.dim}\nolimits}
\def\op{^\mathrm{op}}
\def\Rep{\mathop\mathbf{Rep}\nolimits}
\def\rep{\mathop\mathbf{rep}\nolimits}
\def\Dim{\mathop\mathbf{Dim}\nolimits}
\def\hom{\mathop{\mathcal H\hspace{-.15ex}\mathit{om}}\nolimits}
\def\alg{\mathop\mathrm{Alg}}
\def\Hom{\mathop\mathrm{Hom}}
\def\End{\mathop\mathrm{End}}
\def\Ker{\mathop\mathrm{Ker}}
\def\im{\mathop\mathrm{Im}}
\def\8{\infty}		\def\+{\oplus}
\def\*{\otimes}	\def\dd{\partial}	
\def\id{\mathrm{id}}

\def\al{\alpha}	\def\be{\beta}	\def\ga{\gamma}
\def\Ga{\Gamma}	\def\io{\iota}		\def\la{\lambda}
\def\La{\Lambda}	\def\si{\sigma}	\def\om{\omega}	

\def\fK{\mathbf k}	\def\bA{\mathbf A}	\def\bB{\mathbf B}
\def\fM{\mathbf m}	\def\bI{\mathbf I}	\def\bJ{\mathbf J}
\def\bE{\mathbf E}	\def\fR{\mathbf r}	\def\bR{\mathbf R}
\def\fA{\mathbf a}	\def\fB{\mathbf b}	\def\fV{\mathbf v}
\def\fD{\mathbf d}	\def\bG{\mathbf G}	\def\bK{\mathbf K}	

\def\mN{\mathbb N}	\def\mZ{\mathbb Z}		\def\mA{\mathbb A}
\def\mP{\mathbb P}	\def\mD{\mathbb D}	

\def\kA{\mathcal A}		\def\kB{\mathcal B}
\def\kC{\mathcal C}	\def\kD{\mathcal D}
\def\kH{\mathcal H}	\def\kO{\mathcal O}
\def\kV{\mathcal V}		\def\kW{\mathcal W}
\def\kF{\mathcal F}		\def\kG{\mathcal G}
\def\kQ{\mathcal Q}	\def\kP{\mathcal P}
\def\kI{\mathcal I}		\def\kJ{\mathcal J}
\def\kS{\mathcal S}		\def\kM{\mathcal M}
\def\kL{\mathcal L}		\def\kT{\mathcal T}
\def\kU{\mathfrak U}	\def\rH{\mathrm H}

\def\gM{\mathfrak m}	\def\dA{\mathfrak A}
\def\dB{\mathfrak B}	\def\dC{\mathfrak C}
\def\dD{\mathfrak D}	\def\dS{\mathfrak S}

\def\sV{\mathsf V}		\def\sS{\mathsf S}

\def\oV{\overline{\mathcal V}}	\def\ro{\varrho}
\def\tW{\widetilde{\mathcal W}}		\def\tmu{\tilde\mu}

\begin{document}

\title{Derived tame and derived wild algebras}
 \author{Yuriy A. Drozd}
 \address{Department of Mechanics and Mathematics,
 Kyiv Taras Shevchenko University, 01033 Kyiv, Ukraine}
 \email{yuriy@drozd.org}
 \urladdr{http://drozd.org/\~{}yuriy}

  \begin{abstract}
 We prove that every finite dimensional algebra over an algebraically closed field is either derived tame or derived wild.
 We also prove that any deformation of a derived wild algebra is derived wild.
 \end{abstract} 

\maketitle
\tableofcontents

 \section*{Introduction}

 This is a talk given by the author at the IV International Algebraic Conference in Ukraine (Lviv, August 2003).
 It is devoted to the notions of derived tameness and wildness of algebras and the recent progress in the topic.
 The main result, obtained by Viktor Bekkert and the author, is the tame-wild dichotomy for derived categories,
 which is an analogue of the author's theorem for algebras \cite{d1}. The proof is also very much alike that of
 \cite{d1}; it relies on the technique of matrix problems (representations of boxes) and a reduction algorithm
 for matrices. The principal distinction is that here we have to consider non-free (even non-semi-free) boxes.
 Fortunately, the required class of boxes can be dealt with in a similar way; it is considered in Section \ref{s3}
 of this paper. Section \ref{s1} is devoted to the general notions related to derived categories, derived tameness,
 etc., while Section \ref{s2} presents a transition to matrix problems. Further, we consider the results related to
 deformations and degenerations of derived tame and derived wild algebras obtained by the author \cite{d2}.
 Namely, in Section \ref{s4} we construct some ``almost versal'' families of complexes with projective bases, 
 analogous to the families of modules constructed in \cite{dg}. It makes possible to introduce the ``number of
 parameters'' defining complexes of given rank and to prove (in Section \ref{s5}) that this number is upper
 semi-continuous in flat families of algebras. As a corollary, we get that a deformation of a derived tame algebra
 is derived tame; respectively, a degeneration of a derived wild algebra is derived wild. We also explain the
 situation that can arise if one considers families that are non-flat (especially the example of Br\"ustle \cite{br}).

 Since it is a survey, we sometimes only sketch proofs referring to the papers \cite{bd,d2}, though we try to
 give all necessary definitions, especially related to derived categories and boxes.

 I am grateful to Viktor Bekkert and Igor Burban for fruitful collaboration and useful discussions,
 which were of great influence in preparing this paper. I am also thankful to Birge Huisgen-Zimmermann who
 inspired me to start these investigations.

 \section{Derived categories}
 \label{s1}

 In what follows we consider finite dimensional algebras over an algebraically closed field $\fK$. Let $\bA$
 be such an algebra. As usually, we write $\*$ and $\Hom$ instead of $\*_\fK$ and $\Hom_\fK$; we
 denote by $V^*$ the dual vector space $\Hom(V,\fK)$. All considered categories $\kA$ are also supposed 
 to be \emph{linear categories} over the field $\fK$, which means that all sets $\kA(a,b)$ are vector spaces
 over $\fK$ and the multiplication of morphisms is $\fK$-bilinear.
 Recall that the \emph{derived category} $\kD(\bA)$ of (finite dimensional)
 $\bA$-modules is defined as follows \cite{gm,hap}.
 First consider the category of \emph{complexes} $\kC(\bA)$.
 Its objects are sequences $M_\bp=(M_n,d_n)$ of finite dimensional modules and their homomorphisms
 \begin{equation}\label{e11}
  \begin{CD}
  \dots  @>d_{n+2}>> M_{n+1} @>d_{n+1}>> M_n @>d_n>> M_{n-1} @>d_{n-1}>> \dots  
 \end{CD} 
 \end{equation} 
 such that $d_{n+1}d_n=0$ for all $n$. A \emph{morphism} $\phi_\bp=(\phi_n)$
 of a complex \eqref{e11} to another complex $M'_\bp=(M'_n,d'_n)$ is a commutative diagram
 \begin{equation}\label{e12}
  \begin{CD}
  \dots  @>d_{n+2}>> M_{n+1} @>d_{n+1}>> M_n @>d_n>> M_{n-1} @>d_{n-1}>> \dots  \\
  &&	@V\phi_{n+1}VV  @V\phi_nVV  @V\phi_{n-1}VV \\
  \dots  @>d'_{n+2}>> M'_{n+1} @>d'_{n+1}>> M'_n @>d'_n>> M'_{n-1} @>d'_{n-1}>> \dots  
 \end{CD} 
 \end{equation} 
 One says that such a morphism is \emph{homotopic to zero} if there are homomorphisms
 $\si_n:M_n\to M'_{n+1} \ (n\in\mN)$ such that $\phi_n=\si_{n-1}d_n+d'_{n+1}\si_n$ for all
 $n$. The factor category $\kH(\bA)$ of $\kC(\bA)$ modulo the ideal of morphisms homotopic to zero
 is called the \emph{homotopic category} of $\bA$-modules. For each $n$ the \emph{$n$-th homology}
 of a complex \eqref{e11} is defined as $\rH_n(M_\bp)=\Ker d_n/\im d_{n+1}$. Obviously,
 a morphism $\phi_\bp$ of complexes induces homomorphisms of homologies $\rH_n(\phi_\bp):
 \rH_n(M_\bp)\to\rH_n(M'_\bp)$, and if $\phi_\bp$ is homotopic to zero, it induces
 zero homomorphisms of homologies. Hence, $\rH_n$ can be considered as functors $\kH(\bA)\to\vec$,
 the category of (finite dimensional) vector spaces. One call a morphism $\phi_\bp$ from $\kC(\bA)$
 or from $\kH(\bA)$ \emph{quasi-isomorphism} if all $\rH_n(\phi_\bp)$ are isomorphisms.
 Now the derived category $\kD(\bA)$ is defined as the category of fractions $\kH(\bA)[\kQ^{-1}]$,
 where $\kQ$ is the set of all quasi-isomorphisms.

 One calls a complex \eqref{e11} \emph{right bounded} (\emph{left bounded}, \emph{bounded})
 if there is $n_0$ such that $M_n=0$ for $n<n_0$ (respectively, there is $n_1$ such that $M_n=0$ for
 $n>n_1$, or there are both). The corresponding categories are denoted by $\kC^-,\kH^-,
 \kD^-$ (respectively, by $\kC^+,\kH^+,\kD^+$, or by $\kC^b,\kH^b,\kD^b$). In this paper we 
 mainly deal with the bounded derived category $\kD^b(\bA)$.
 The category $\bA\md$ of (finite dimensional) $\bA$-modules naturally embeds into $\kD(\bA)$ (even
 in $\kD^b(\bA)$): a module $M$ is identified with the complex $M_\bp$ such that $M_0=M,\,
 M_n=0$ for $n\ne0$. 

 Since the category of modules has enough projectives, one can replace, when considering right bounded
 homotopic or derived category, arbitrary complexes by projective ones, i.e. consisting of projective modules
 only. We denote by $\kP^-(\bA)$ and by $\kP^b(\bA)$ the categories of right bounded and bounded
 projective complexes. Actually, $\kD^-(\bA)\simeq\kH^-(\bA)\simeq\kP^-(\bA)/I_h$, where $I_h$ is the
 ideal of morphisms homotopic to zero \cite{gm,hap}. Moreover, in finite dimensional case every module
 $M$ has a \emph{projective cover}, i.e. there is an epimorphism $p(M):P(M)\to M$, where $P(M)$
 is projective, such that $\Ker p(M)\subseteq\rad P(M)$, the radical of $P(M)$ \cite{dk}. Therefore, one
 can only consider \emph{minimal complexes}, i.e. projective complexes such that $\im d_n\subseteq\rad
 P_{n-1}$ for all $n$. We denote by $\kP^-_{\min}(\bA)$ the category of minimal projective complexes.
 Thus $\kD^-(\bA)\simeq\kP^-_{\min}(\bA)/I_h$. One immediately checks that the image in $\kD^-(\bA)$ of a
 morphism $\phi_\bp$ of  minimal complexes is an isomorphism \iff $\phi$ itself is an isomorphism.
 Hence, if we are interested in classification of objects from $\kD^-(\bA)$, we can replace it by
 $\kP^-_{\min}(\bA)$. Unfortunately, it is no more the case if we consider $\kD^b(\bA)$. Certainly,
 if $\gdim\bA<\8$, one has that $\kD^b(\bA)\simeq\kP^b_{\min}(\bA)/I_h$, since any bounded
 complex has a bounded projective resolution. On the contrary, if $\gdim\bA=\8$, it is wrong even for
 modules. Nevertheless, we propose the following approximation of $\kD^b(\bA)$, which is enough 
 for classification purpose. (Compare it with \cite{hs}.)

 We denote by $\kP^m(\bA)$ the category of bounded projective complexes $M_\bp$ such that 
 $M_n=0$ for $n>m$ (note that the left bound is not prescribed). We say that a morphism $\phi_\bp$
 in $\kP^m(\bA)$ is \emph{quasi-homotopic to zero} if there are homomorphisms
 $\si_k:M_n\to M'_{n+1}$ such that $\phi_n=\si_{n-1}d_n+d'_{n+1}\si_n$ for all $n<m$. Let
 $\kH^m(\bA)=\kP^m(\bA)/I_{qh}$, where $I_{qh}$ is the ideal of morphisms quasi-homotopic to
 zero. Especially, if $P_\bp$ is a minimal complex from $\kP^m(\bA)$, it is isomorphic in $\kH^m(\bA)$
 to a minimal complex $\tilde P_\bp$ such that $\tilde P_n=P_n,\,\tilde d_n=d_n$ for $n\ne m$, while $P_m=
 \tilde P_m\+Q$ so that $Q\subseteq\Ker d_m$, $\tilde d_m=d_m|_{\tilde P_m}$ and $\Ker \tilde d_m\subseteq\rad \tilde P_m$.
  There is a natural functor $\kH^m(\bA)\to\kH^{m+1}(\bA)$, which maps
 a complex $P_\bp$ to the complex $P^+_\bp$, where $P^+_n=P_n,\ d^+_n=d_n$ for
 $n\le m$, while $P^+_{m+1}=P(\Ker d_m)$ and $d^+_{m+1}$ is the epimorphism $p(\Ker d_m):
 P^+_{m+1}\to\Ker d_m$. Conversely, there is a functor $\kH^{m+1}(\bA)\to\kH^m(\bA)$, which
 just abandon $P_{m+1}$ in a complex $P_\bp$. One easily verifies the following results. 

  \begin{prop}\label{11} 
  \begin{enumerate} 
 \item 	$\kD^b(\bA)\simeq\varinjlim_m\kH^m(\bA)$; 
 \item	$\kD^-(\bA)\simeq\varprojlim_m\kH^m(\bA)$. 
 \end{enumerate} 
 Moreover, if complexes $P_\bp$ and $P'_\bp\in\kP^m(\bA)$ are minimal, their
 images in $\kD^b(\bA)$ are isomorphic \iff $\tilde P_\bp\simeq\tilde P'_\bp$ as complexes. 
 \end{prop} 

 Thus dealing with classification problems, one can deal with $\kH^m$ or with $\kP^m$  instead of $\kD^b$.
 Note also that if we fix all $\rH_n(P_\bp)$ for a complex from $\kP^m_{\min}(\bA)$, there are finitely many
 possibilities for the values of $P_n$. 

 Taking into consideration these remarks, one can define \emph{derived wild} and \emph{derived tame} algebras
 as follows. 

  \begin{defin}\label{wild}
  \begin{enumerate}
\item 	Let $\bR$ be a $\fK$-algebra. A \emph{family of $\bA$-complexes} based on $\bR$
 is a complex of finitely generated projective $\bA\*\bR\op$-modules $P_\bp$. We denote by
 $\kP^m(\bA,\bR)$ the category of all bounded families with $P_n=0$ for $n>m$ (again we do not prescribe
 the right bound). For such a family $P_\bp$ and an $\bR$-module $L$ we denote by $P_\bp(L)$ the
 complex $(P_n\*_\bR L,d_n\*1)$. If $L$ is finite dimensional, $P_\bp(L)\in\kP^m(\bA)$.

 \item	We call a family $P_\bp$ \emph{strict} if for every finite dimensional $\bR$-modules $L,L'$
  \begin{enumerate}
  \item   $P_\bp(L)\simeq P_\bp(L')$ \iff $L\simeq L'$;
  \item	  $P_\bp(L)$ is indecomposable \iff so is $L$.
  \end{enumerate}

 \item	We call $\bA$ \emph{derived wild} if it has a strict family of complexes over every finitely generated
 $\fK$-algebra $\bR$.
 \end{enumerate}
 \end{defin} 

 The following fact is well known.

  \begin{prop}\label{test}
  An algebra $\bA$ is derived wild \iff it has a strict family over one of the following algebras:
 \begin{enumerate}
\item 	free algebra $\fK\langle x,y\rangle$ in two variables;
 \item	polynomial algebra $\fK[x,y]$ in two variables;
 \item	power series algebra $\fK[[x,y]]$ in two variables.
\end{enumerate}
 \end{prop} 

   Let $\lst As$ be a set of representatives of isomorphism classes of indecomposable projective $\bA$-modules.
 (For instance, one can choose a decomposition $1=\sum_{i=1}^me_i$, where $e_i$ are primitive orthogonal idempotents,
 take $\lst es$ all pairwise non-conjugate of them, and set $A_i=\bA e_i$.) If $P$ is an arbitrary projective $\bA$-module,
 it uniquely decomposes as $P\simeq\bigoplus_{i=1}^sr_iA_i$. Set $\fR(P)=\row rs$ and call it the \emph{vector rank} of
 $P$. For a projective complex $P_\bp$ set $\fR_\bp(P_\bp)=(\fR(P_n)\,|\,n\in\mZ)$. Given an arbitrary sequence
 $\fV_\bp=(\fV_n\,|\,n\in\mZ)$ denote by $\kP(\fV_\bp,\bA)$ the set of complexes $P_\bp$ such that $\fR_\bp(P_\bp)=\fV_\bp$. 

  \begin{defin}\label{tame}
  \begin{enumerate}
 \item   A \emph{rational algebra} is a $\fK$-algebra $\fK[t,f(t)^{-1}]$ for a non-zero polynomial $f(t)$.
 A \emph{rational family} of $\bA$-complexes is a family over a rational algebra $\bR$.
 \item	For a rational family $P_\bp$ denote by $\fR_\bp(P_\bp)=\fR_\bp(P_\bp(L))$ , where $L$ is a one-dimensional
 $\bR$-module. (This value does not depend on the choice of $L$.)
 \item	An algebra $\bA$ is called \emph{derived tame} if there is a set of rational families of bounded
 $\bA$-complexes $\dS$ such that:
	\begin{enumerate}
	\item 	for every (bounded) $\fV_\bp$ the set $\dS(\fV_\bp)
	=\setsuch{P_\bp\in\dS}{\fR_\bp(P_\bp)=\fV_\bp}$ is finite.
	\item	for every $\fV_\bp$ all indecomposable complexes from $\kP(\fV_\bp,\bA)$, except finitely many
 	of them (up to isomorphism) are isomorphic to a complex $P_\bp(L)$ for some $P_\bp\in\dS$ 	
	and some finite dimensional $L$.
	\end{enumerate}
  \end{enumerate} 
 We call $\dS$ a \emph{parameterising set} of $\bA$-complexes. 
 \end{defin} 
 
 These definitions do not formally coincide with other definitions of derived tame and derived wild algebras,
 for instance, those proposed in \cite{ge,gk}, but all of them are evidently equivalent.
 It is obvious (and easy to prove) that neither algebra can be both derived tame and derived wild. 
 The following result (``tame-wild dichotomy for derived categories'')
 has recently been proved by V.\,Bekkert and the author \cite{bd}.

  \begin{theorem}[\sc Main Theorem]\label{main}
 Every finite dimensional algebra over an algebraically closed field is either derived tame or derived wild.  
 \end{theorem} 

 \section{Reduction to matrix problems}
 \label{s2}

 We recall now the main notions related to the matrix problems (representations of boxes). More detailed
 exposition can be found in \cite{ds}. A \emph{box} is a pair
 $\dA=(\kA,\kV)$, where $\kA$ is a ($\fK$-linear) category and $\kV$ is an $\kA$-coalgebra, i.e. an $\kA$-bimodule
 supplied with \emph{comultiplication} $\mu:\kV\to\kV\*_\kA\kV$ and \emph{counit} $\io:\kV\to\kA$, which are
 homomorphisms of $\kA$-bimodules and satisfy the usual coalgebra conditions  
 $$ 
   (\mu\*1)\mu=(1\*\mu)\mu,\quad  i_l(\io\*1)\mu=i_r(1\*\io)\mu=\id,
 $$ 
 where $i_l:\kA\*_\kA\kV\simeq\kV$ and $i_r:\kV\*_\kA\kA\simeq\kV$ are the natural isomorphisms. The kernel $\oV=\Ker\io$
 is called the \emph{kernel of the box}. A \emph{representation}
 of such a box in a category $\kC$ is a functor $M:\kA\to\kC$. Given another representation $N:\kA\to\kC$, a \emph{morphism}
 $f:M\to N$ is defined as a homomorphism of $\kA$-modules $\kV\*_\kA M\to N$, or, equivalently, as a homomorphism of
 $\kA$-bimodules $\kV\to\Hom_\kC(M,N)$, the latter supplied with the obvious $\kA$-bimodule structure. The composition
 $gf$ of $f:M\to N$ and $g:N\to L$ is defined as the composition 
 $$ 
   \begin{CD}
  \kV\*_\kA M@>\mu\*1>> \kV\*_\kA\kV\*_\kA M @>1\*f>>\kV\*_\kA N @>g>> L,
 \end{CD} 
 $$ 
 while the identity morphism $\id_M$ of $M$ is the composition 
 $$ 
    \begin{CD}
  \kV\*_\kA M@>\io\*1>> \kA\*_\kA M @>i_l >> M.
 \end{CD} 
 $$ 
 Thus we obtain the \emph{category of representations} $\Rep(\dA,\kC)$. If $\kC=\vec$, we just write
 $\Rep(\dA)$. If $f$ is a morphism and $\ga\in\kV(a,b)$, we denote by $f(\ga)$ the morphism
 $f(b)(\ga\*\_\,): M(a)\to N(a)$.
 A box $\dA$ is called \emph{normal} (or \emph{group-like}) if there is a set of elements
 $\om=\setsuch{\om_a\in\kV(a,a)}{a\in\ob\kA}$ such that $\io(\om_a)=1_a$ and $\mu(\om_a)=
 \om_a\*\om_a$ for every $a\in\ob\kA$. In this case, if $f$ is an isomorphism, all morphisms $f(\om_a)$
 are isomorphisms $M(a)\simeq N(a)$. This set is called a \emph{section} of $\dA$. For a normal box, one defines
 the \emph{differentials} $\dd_0:\kA\to\oV$ and $\dd_1:\oV\to\oV\*_\kA\oV$ setting
 \begin{align*} 
  \dd_0(\al)&=\al\om_a-\om_b\al \quad \text{ for } \, \al\in\kA(a,b);\\
  \dd_1(\ga)&=\mu(\ga)-\ga\*\om_a-\om_b\*\ga \quad \text{ for }\, \ga\in\oV(a,b). 
 \end{align*} 
 Usually we omit indices, writing $\dd\al$ and $\dd\ga$. 

 Recall that a \emph{free category} $\fK\Ga$, where $\Ga$ is an oriented graph, has the vertices of $\Ga$
 as its objects  and the paths from $a$ to $b$ ($a,b$ being two vertices) as a basis of the vector space
 $\fK\Ga(a,b)$. A \emph{semi-free} category is a category of fractions $\fK\Ga[S^{-1}]$, where
 $S=\setsuch{g_\al(\al)}{\al \text{ runs through a set of loops}}$
  (called the \emph{marked loops}) of the graph $\Ga$. The arrows of $\Ga$ are called the
 \emph{free} (respectively, \emph{semi-free}) generators of the free (semi-free) category. A \emph{normal box}
 $\dA=(\kA,\kV)$ is called \emph{free} (\emph{semi-free}) if such is the category $\kA$, moreover,
 the kernel $\oV=\Ker\io$ of the box is a free $\kA$ -bimodule and $\dd\al=0$
 for each marked loop $\al$. A set of \emph{free} (respectively, \emph{semi-free}) {generators} of such a
 box is a union $\sS=\sS_0\cup\sS_1$, where $\sS_0$ is a set of free (semi-free) generators of the
 category $\kA$ and $\sS_1$ is a set of free generators of the $\kA$-bimodule $\kV$. A set of
 free (or semi-free) generators $\sS$ is called \emph{triangular} if there is a function $\ro:\sS\to\mN$
 (we call it ``\emph{weight}'')
 such that, for every $\al\in\sS$, the differential $\dd\al$ belongs to the sub-box generated by
 $\setsuch{\be\in\sS}{\ro(\be)<\ro(\al)}$; especially, if $\ro(\al)=0$, $\al$ is \emph{minimal}, i.e.
 $\dd\al=0$. A free (semi-free) box having a triangular set of free (semi-free) generators is called \emph{triangular}.
  For a triangular box a morphism of representations $f:M\to N$ is an isomorphism \iff all $f(\om_a)$ are isomorphisms. 
 In what follows we always suppose that our graphs are \emph{locally finite}, i.e. only have finitely many arrows
 starting or ending at a given vertex. If such a graph $\Ga$ has no oriented cycles, then the category $\fK\Ga$ is
 \emph{locally finite dimensional}, i.e. all its spaces of morphisms are finite dimensional.

 We call a category $\kA$ \emph{trivial} if it is a free category generated by a trivial graph (i.e. one with no arrows);
 thus $\kA(a,b)=0$ if $a\ne b$ and $\kA(a,a)=\fK$. We call $\kA$ \emph{minimal}, if it is a semi-free category
 with a set of semi-free generators consisting of loops only, at most one loop at each vertex. Thus $\kA(a,b)=0$ again
 if $a\ne b$, while $\kA(a,a)$ is either $\fK$ or a rational algebra. We call a normal box $\dA=(\kA,\kV)$
 \emph{so-trivial} if $\kA$ is trivial, and \emph{so-minimal} if $\kA$ is minimal and all its loops $\al$
 are minimal too (i.e. with $\dd\al=0$). 
 
 In \cite{d1} (cf. also \cite{ds}) the classification of representations of an arbitrary finite dimensional algebra was reduced
 to representations of a free triangular box. To deal with derived categories we have to consider a wider class of boxes.
 First, a \emph{factor-box}  of a box $\dA=(\kA,\kV)$ modulo an ideal $\kI\subseteq\kA$ is defined as
 the box $\dA/\kI=(\kA/\kI,\kV/(\kI\kV+\kV\kI)$ (with obvious comultiplication and counit). Note that if
 $\dA$ is normal, so is $\dA/\kI$. 

  \begin{defin}\label{slice}
 A \emph{sliced box} is a factor-box $\dA/\kI$, where $\dA=(\kA,\kV)$ is a free box such that the set of its objects
 $\sV=\ob\kA$ is a disjoint union $\sV=\bigcup_{i\in\mZ}\sV_i$ such that $\kA(a,b)=0$ if $a\in\sV_i,\,b\in\sV_j$ with
 $j>i$, $\kA(a,a)=\fK$, and $\kV(a,b)=0$ if $a\in\sV_i,\,b\in\sV_j$ with $i\ne j$.  We call a sliced box
 $\dA/\kI$ \emph{triangular} if so is the free box $\dA$. The partition $\sV=\bigcup_i\sV_i$ is called a \emph{slicing}.
 \end{defin} 

 Certainly, in this definition we may assume that the elements of the ideal $\kI$ are linear combinations of paths of
 length at least $2$. Otherwise such an element is a linear combination of arrows, so we can just eliminate one
 of these arrows from the underlying graph.

 Note that for every representation $M\in\Rep(\dA)$, where $\dA$ is a free (semi-free, sliced) box with the set of objects
 $\sV$, one can consider its \emph{dimension} $\Dim(M)$, which is a function $\sV\to\mN$, namely
 $\Dim(M)(a)=\dim M(a)$. We call such a representation \emph{finite dimensional} if its support
 $\supp M=\setsuch{a\in\sV}{M(a)\ne0}$ is finite and denote by $\rep(\dA)$ the category of finite
 dimensional representations. Having these notions, one can easily reproduce the definitions of families of
 representations, especially strict families, wild and tame boxes; see \cite{d1,ds} for details. The following
 procedure, mostly modelling that of \cite{d1}, let us replace derived categories by representations of
 sliced boxes.

 Let $\bA$ be a finite dimensional algebra, $\bJ$ be its radical. As far as we are interested in $\bA$-modules
 and complexes, we can replace $\bA$ by a Morita equivalent reduced algebra, thus suppose that $\bA/\bJ\simeq
 \fK^s$ \cite{dk}. Let $1=\sum_{i=1}^s$, where $e_i$ are primitive orthogonal idempotents; set
 $\bA_{ji}=e_j\bA e_i$ and $\bJ_{ji}=e_j\bJ e_i$; note that $\bJ_{ji}=\bA_{ji}$ if $i\ne j$.
 We denote by $\kS$ the trivial category with the set of objects
 $\setsuch{(i,n)}{n\in\mN,\,i=1,2,\dots,s}$ and consider the $\kS$-bimodule $\kJ$ such that  
 $$ 
   \kJ\big((i,n),(j,m)\big)= \begin{cases}
  0 &\text{if } m\ne n-1,\\
  \bJ_{ji}^* &\text{if }m=n-1.
 \end{cases} 
 $$ 
 Let $\kB=\kS[\kJ]$ be the tensor category of this bimodule; equivalently, it is  the free category having the same
 set of objects as $\kS$ and the union of bases of all $\kJ\big((i,n),(j,m)\big)$ as a set of free generators. 
 Denote by $\kU$ the $\kS$-bimodule such that 
  $$
 \kU\big((i,n),(j,m)\big)= \begin{cases}
  0 &\text{if } n\ne m,\\
 \bA_{ji}^* &\text{if } n=m
 \end{cases} 
 $$
 and set $\tW=\kB\*_\kS\kU\*_\kS\kB$. Dualizing the multiplication $\bA_{kj}\*\bA_{ji}\to\bA_{ki}$, we
 get homomorphisms
 \begin{align*}
  \la_r:& \kB\to \kB\*_\kS\tW,\\
  \la_l:& \kB\to \tW\*_\kS\kB,\\
  \tmu:& \tW\to \tW\*_\kS\tW.
 \end{align*} 
 In particular, $\tmu$ defines on $\tW$ a structure of $\kB$-coalgebra. Moreover, the sub-bimodule $\kW_0$
 generated by $\im(\la_r-\la_l)$ is a coideal in $\tW$, i.e. $\tmu(\kW_0)\subseteq\kW_0\*_\kB\tW\+\tW\*_\kB\kW_0$.
 Therefore, $\kW=\tW/\kW_0$ is also a $\kB$-coalgebra, so we get a box $\dB=(\kB,\kW)$. One easily checks
 that it is free and triangular. 

 Dualizing multiplication also gives a mapping
 \begin{equation}\label{e21} 
 \nu:\bJ_{ji}^*\to\bigoplus_{k=1}^s\bJ_{jk}^*\*\bJ_{ki}^*. 
 \end{equation} 
 Namely, if we choose bases $\set\al,\,\set\be\,\set\ga$ in the spaces, respectively, $\bJ_{ji},$ $\bJ_{jk},\,\bJ_{ki}$,
 and dual bases $\set{\al^*},\,\set{\be^*},\,\set{\ga^*}$ in their duals, then $\be^*\*\ga^*$ occurs in $\nu(\al^*)$
 with the same coefficient as $\al$ occurs in $\be\ga$.
 Note that the right-hand space in \eqref{e21} coincide with each $\kB\big((i,n),(j,n-2)\big)$. Let $\kI$ be the ideal
 in $\kB$ generated by the images of $\nu$ in all these spaces and $\dD=\dB/\kI=(\kA,\kV)$, where
 $\kA=\kB/\kI,\ \kV=\kW/(\kI\kW+\kW\kI)$. If necessary, we write $\dD(\bA)$ to emphasise that this box has been 
 constructed from a given algebra $\bA$. Certainly, $\dD$ is a sliced triangular box, and the following result holds.

  \begin{theorem}\label{box}
  The category of finite dimensional representations $\rep(\dD(\bA))$ is equivalent to the category $\kP^b_{\min}(\bA)$
 of bounded minimal projective $\bA$-complexes.
 \end{theorem} 
 \begin{proof}
 Let $A_i=\bA e_i$; they form a complete list of non-isomorphic indecomposable projective $\bA$-modules; set also
 $J_i=\rad A_i=\bJ e_i$. Then $\Hom_\bA(A_i,J_j)\simeq\bJ_{ji}$. A representation $M\in\rep(\dD)$ is given by
 vector spaces $M(i,n)$ and linear mappings 
 $$ 
 M_{ji}(n):\bJ_{ji}^*=\kA\big((i,n),(j,n-1)\big)\to\Hom\big(M(i,n),M(j,n-1)\big)
 $$
 subject to the relations 
 \begin{equation}\label{e22}
   \sum_{k=1}^s \fM\big(M_{jk}(n)\*M_{ki}(n+1)\big)\nu(\al)=0
 \end{equation}
 for all $i,j,k,n$ and all $\al\in\bJ_{ji}$, where $\fM$ denotes the multiplication of mappings 
  
 \begin{multline*} 
   \Hom\big(M(k,n),M(j,n-1)\big)\*\Hom\big(M(i,n+1),M(k,n)\big)\to\\
	\to\Hom\big(M(i,n+1),M(j,n-1)\big).   
 \end{multline*}  
 For such a representation, set $P_n=\bigoplus_{i=1}^s A_i\*M(i,n)$. Then $\rad P_n=\bigoplus_{i=1}^n J_i\*M(i,n)$
 and
 \begin{align*}
  \Hom_\bA(P_n,\rad P_{n-1})&\simeq \bigoplus_{i,j} \Hom_\bA\big(A_i\*M(i,n),J_j\*M(j,n-1)\big)\simeq\\
	&\simeq \bigoplus_{ij} \Hom\big(M(i,n),\Hom_\bA\big(A_i,J_j\*M(j,n-1)\big)\big)\simeq\\
	&\simeq \bigoplus_{ij} M(i,n)^*\*\bJ_{ji}\*M(j,n-1) \simeq\\
	&\simeq \bigoplus_{ij} \Hom\big(\bJ^*_{ji},\Hom\big(M(i,n),M(j,n-1)\big)\big).
 \end{align*} 
 Thus the set $\setsuch{M_{ji}(n)}{i,j=1,2,\dots,s}$ defines a homomorphism $d_n:P_n\to P_{n-1}$ and vice versa.
 Moreover, one easily verifies that the condition \eqref{e22} is equivalent to the relation $d_nd_{n+1}=0$. Since
 every projective $\bA$-module can be given in the form $\bigoplus_{i=1}^sA_i\*V_i$ for some uniquely defined
 vector spaces $V_i$, we get a \oc between finite dimensional representations of $\dD$ and bounded minimal
 complexes of projective $\bA$-modules. In the same way one also establishes \oc between morphisms of
 representations and of the corresponding complexes, compatible with their multiplication,
 which accomplishes the proof.
  \end{proof} 

  \begin{corol}\label{twbox}
  An algebra $\bA$ is derived tame (derived wild) if so is the box $\dD(\bA)$. 
 \end{corol} 

 \section{Proof of the main theorem}
 \label{s3}

 Now we are able to prove the main theorem. Namely, according to Corollary \ref{twbox}, it follows from the
 analogous result for sliced boxes.

  \begin{theorem}\label{mbox}
 Every sliced triangular box is either tame or wild.   
 \end{theorem} 

 Actually, just as in \cite{d1} (see also \cite{ds}), we shall prove this theorem in the following form.

\theoremstyle{plain}
\newtheorem*{31a}{Theorem 3.1a}
  \begin{31a}\label{mini}
  Suppose that a sliced triangular box $\dA=(\kA,\kV)$ is not wild. For every dimension $\fD$ of its
 representations there is a functor $F_\fD:\kA\to\kM$, where $\kM$ is a minimal category, such that
 every representation $M:\kA\to\vec$ of $\dA$ of dimension $\Dim (M)\le\fD$ is isomorphic to the
 inverse image $F^*N=N\circ F$ for some  functor $N:\kM\to\vec$. Moreover, $F$ can be chosen
 \emph{strict}, which means that $F^*N\simeq F^*N'$ implies $N\simeq N'$ and $F^*N$ is
 indecomposable if so is $N$.
 \end{31a} 

  \begin{erem}\label{r21}
  We can consider the induced box $\dA^F=(\kM,\kM\*_\kA\kV\*_\kA\kM)$. It is a so-minimal box, and $F^*$
 defines a full and faithful functor $\rep(\dA^F)\to\rep(\dA)$. Its image consists of all representations 
 $M:\kA\to\vec$ that factorise through $F$. 
 \end{erem} 

 \begin{proof}
 As we only consider finite dimensional representations, we may assume that the set of
 objects is finite. Hence the slicing $\sV=\bigcup_i\sV_i$ (see Definition \ref{slice}) is finite too:
 $\sV=\bigcup_{i=1}^m\sV_i$ and we use induction by $m$. If $m=1$, $\dA$ is free, and
 our claim has been proved in \cite{d1}. So we may suppose that the theorem is true for smaller
 values of $m$, especially, it is true for the restriction $\dA'=(\kA',\kV')$ of the box $\dA$
 onto the subset $\sV'=\bigcup_{i=2}^m\sV_i$. Thus there is a strict functor $F':\kA'\to\kM$,
 where $\kM$ is a minimal category, such that every representation of $\dA'$ of dimension smaller
 than $\fD$ is of the form ${F'}^*N$ for $N:\kM\to\vec$. Consider now the amalgamation
 $\kB=\kA\bigsqcup^{\kA'}\kM$ and the box $\dB=(\kB,\kW)$, where
 $\kW=\kB\*_\kA\kV\*_\kA\kB$. The functor $F'$ extends to a functor $F:\kA\to\kB$ and
 induces a homomorphism of $\bA$-bimodules $\kV\to\kW$; so it defines a functor $F^*:
 \rep(\dB)\to\rep(\dA)$, which is full and faithful. Moreover, every representation of $\dA$
 of dimension smaller than $\fD$ is isomorphic to $F^*N$ for some $N$, and all possible
 dimensions of such $N$ are restricted by some vector $\fB$.
 Therefore, it is enough to prove the claim of the theorem for the box $\dB$. 

 Note that the category $\kB$ is generated by the loops from $\kM$ and the images of arrows from
 $\kA(a,b)$ with $b\in\sV_1$ (we call them \emph{new arrows}). It implies that all possible
 relations between these morphisms are of the form $\sum_\be \be g_\be(\al)$, where
 $\al\in\kB(a,a)$ is a loop (necessarily minimal, i.e. with $\dd\al=0$), $g_\be$ are
 some polynomials, and $\be$ runs through
 the set of new arrows from $a$ to $b$ for some $b\in\sV_1$. Consider all of these relations for
 a fixed $b$; let them be $\sum_\be \be g_{\be,k}(\al)$. Their coefficients form a matrix
 $\big(g_{\be,k}(\al)\big)$. Making linear transformations of the set $\set\be$ and of the set
 of relations, we can make this matrix diagonal, i.e. make all relations being $\be f_\be(\al)=0$
 for some polynomials $f_\be$. If one of $f_\be$ is zero, the box $\dB$ has a sub-box  
$$
 \xymatrix{
	{a} \ar@(ul,dl)[]_{\al} \ar[rr]^{\be} && b	},
$$

 \medskip\noindent
 with $\dd\al=\dd\be=0$,
 which is wild; hence $\dB$ and $\dA$ are also wild. Otherwise, let $f(\al)\ne0$ be a common multiple
 of all $f_\be(\al)$, $\La=\set{\lst \la r}$ be the set of roots of $f(\al)$. If $N\in\rep(\dB)$ is such
 that $N(\al)$ has no eigenvalues from $\La$, then $f(N(\al))$ is invertible; thus $N(\be)=0$ for all
 $\be:a\to b$. So we can apply the \emph{reduction of the loop} $\al$ with respect to the set
 $\La$  and the dimension $d=\fB(a)$, as in \cite[Propositions 3,4]{d1} or \cite[Theorem 6.4]{ds}.
 It gives a new box that has the same number of loops as $\dB$,
 but the loop corresponding to $\al$ is ``isolated,'' i.e. there are no more arrows starting or ending
 at the same vertex. In the same way we are able to isolate all loops, obtaining a semi-free triangular box 
 $\dC$ and a morphism $G:\dB\to\dC$ such that $G^*$ is full and faithful and all representations of
 $\dB$ of dimensions smaller than $\fB$ are of the form $G^*L$. As the theorem is true for semi-free
 boxes, it accomplishes the proof.
 \end{proof} 

  \begin{erem}
  Applying reduction functors, like in the proof above, we can also extend to sliced boxes (thus to derived
 categories) other results obtained before for free boxes. For instance, we mention 
 the following theorem, quite analogous to that of Crawley-Boevey \cite{cb1}.

  \begin{theorem}\label{generic}
  If an algebra $\bA$ is derived tame, then, for any vector $\fD=(d_n\,|\,n\in\mZ)$ such that almost all
 $d_n=0$, there is at most finite set of \emph{generic $\bA$-complexes} of endolength $\fD$, i.e.
 such indecomposable minimal bounded complexes $P_\bp$ of projective $\bA$-modules, not all
 of which are finitely generated, that $\,\mathrm{length}_\bE(P_n)=d_n$ for all $n$,
 where $\bE=\End_\bA(P_\bp)$.
 \end{theorem} 
 Its proof reproduces again that of \cite{cb1}, with obvious changes necessary to include sliced boxes
 into consideration.
 \end{erem} 

 \section{Families of complexes}
 \label{s4}
 
 We consider now \emph{algebraic families} of $\bA$-complexes, i.e. flat families over an algebraic variety $X$.
 Such a family is a complex $\kF_\bp=(\kF_n,d_n)$ of flat coherent $\bA\*\kO_X$-modules. We always assume
 this complex bounded and \emph{minimal}; the latter means that $\im d_n\subseteq\bJ\kF_{n-1}$ for all $n $,
 where $\bJ=\rad\bA$. We also assume that $X$ is connected; it implies that the vector rank
 $\fR_\bp(\kF_\bp(x))$ is constant, so we can call it the \emph{vector rank of the family} $\kF$ and
 denote it by $\fR_\bp(\kF_\bp)$ Here, as usually, $\kF(x)=\kF_x/\gM_x\kF_x$, where $\gM_x$
 is the maximal ideal of the ring $\kO_{X,x}$. We call a family $\kF_\bp$ \emph{non-degenerate}
 if, for every $x\in X$, at least one of  $d_n(x):\kF_n(x)\to\kF_{n-1}(x)$ is non-zero. 
 Having a family $\kF_\bp$ over $X$ and a regular mapping $\phi:Y\to X$, one gets the inverse image $\phi^*(\kF)$,
 which is a family of $\bA$-complexes over the variety $Y$ such that $\phi^*(\kF)(y)\simeq\kF(\phi(y))$.
 If $\kF_\bp$ is non-degenerate, so is $\phi^*(\kF)$. Given an ideal $\bI\subseteq\bJ$, we call a family $\kF_\bp$
 an \emph{$\bI$-family} if $\im d_n\subseteq\bI\kF_{n-1}$ for all $n$. Then any inverse image
 $\phi^*(\kF)$ is an $\bI$-family as well. Just as in \cite{dg}, we construct some ``almost versal''
 non-degenerate $\bI$-families.

 Consider again a complete set of non-isomorphic indecomposable projective $\bA$-modules $\set{\lst As}$.
 For each vector $\fR=\row rs$ set $\fR\bA=\bigoplus_{i=1}^sr_iA_i$, and denote
 $\bI(\fR,\fR')=\Hom_\bA(\fR\bA,\bI\cdot\fR'\bA)$, where $\bI$ is an ideal contained in $\bJ$.
 Fix a vector rank of bounded complexes $\fR_\bp=(\fR_k\,|\,m\le k\le n)$
 and set $\kH=\kH(\fR,\bI)=\bigoplus_{k=m+1}^n\bI(\fR_k,\fR_{k-1})$.
 Consider the projective space $\mP=\mP(\fR_\bp,\bI)=\mP(\kH)$ and its closed subset
 $\mD=\mD(\fR_\bp,\bI)\subseteq\mP$ consisting of all sequences $(h_k)$
 such that $h_{k+1}h_k=0$ for all $k$. Because of the universal property of projective spaces
 \cite[Theorem II.7.1]{ha}, the embedding $\mD(\fR_\bp,\bI) \to\mP(\fR_\bp,\bI)$ gives rise to 
 a non-degenerate $\bI$-family $\kV_\bp=\kV_\bp(\fR_\bp,\bI)$: 
 \begin{equation}\label{e41} 
 \begin{CD} 
    \kV_\bp:\quad \kV_n@>{d_n}>>\kV_{n-1}@>{d_{n-1}}>>\dots\longrightarrow \kV_m ,
 \end{CD} 
 \end{equation} 
 where $\kV_k=\kO_\mD(n-k)\*\fR_k\bA$ for all $m\le k\le n$.
 We call $\kV_\bp(\fR_\bp,\bI)$ the \emph{canonical $\bI$-family of $\bA$-complexes} over $\mD(\fR_\bp,\bI)$.
 Moreover, regular mappings $\phi:X\to\mD(\fR_\bp,\bI)$ correspond to non-degenerate $\bI$-families $\kF_\bp$ with
 $\kF_k=0$ for $k>n$ or $k<m$ and $\kF_k=\kL^{\*(n-k)}\*\fR_k\bA$ for some invertible sheaf $\kL$
 over $X$. Namely, such a family can be obtained as $\phi^*(\kV_\bp)$ for a uniquely defined regular mapping
 $\phi$. Moreover, the following result holds, which shows the ``almost versality'' of the families $\kV_\bp(\fR_\bp,\bI)$.

  \begin{prop}\label{vers}
  For every non-degenerate family of $\bI$-complexes $\kF_\bp$ of vector rank $\fR_\bp$ over an algebraic
 variety $X$, there is a finite open covering $X=\bigcup_jU_j$ such that the restriction of $\kF_\bp$ onto each $U_j$
 is isomorphic to $\phi_j^*\kV_\bp(\fR_\bp,\bI)$ for a regular mapping $\phi_j:U_j\to \mD(\fR_\bp,\bI)$.
 \end{prop}
  \begin{proof}
 For each $x\in X$  there is an open neighbourhood $U\ni x$ such that all restrictions $\kF_k|_U$ are isomorphic to
 $\kO_U\*\fR_k\bA$; so the restriction $\kF_\bp|_U$ is obtained from a regular mapping $U\to\mD(\fR_\bp,\bI)$. 
 Evidently it implies the assertion.
 \end{proof} 

 Note that the mappings $\phi_j$ are not canonical, so we cannot glue them into a ``global'' mapping
 $X\to\mD(\fR_\bp,\bI)$.

 Consider now the group $\bG=\bG(\fR_\bp)=\prod_k\Aut(\fR_k\bA)$, which acts on $\kH(\fR_\bp,\bI)$:
 $(g_k)\cdot(h_k)=(g_{k-1}h_kg_k^{-1})$. It induces the action of $\bG(\fR_\bp)$ on $\mP(\bR_\bp,\bI)$
 and on $\mD(\fR_\bp,\bI)$. The definitions immediately imply that
 $\kV_\bp(\fR_\bp,\bI)(x)\simeq\kV_\bp(\fR_\bp,\bI)(x')$ ($x,x'\in\mD)$ \iff $x$ and $x'$ belong to the same
 orbit of $\bG$. Consider the sets 
 $$
 \mD_i=\mD_i(\fR_\bp,\bI)=\setsuch{x\in\mD}{\dim\bG x\le i}. 
 $$ 
 It is known that they are closed (it follows from the theorem on dimensions of fibres, cf.
 \cite[Exercise II.3.22]{ha} or \cite[Ch.\,I,\,\S\,6,\,Theorem 7]{sh}). We set 
 $$ 
   \pam(\fR_\bp,\bI,\bA)=\max_i\set{\dim\mD_i(\fR_\bp,\bI)-i}
 $$ 
 and call this integer the \emph{parameter number} of $\bI$-complexes of vector rank $\fR_\bp$.
 Obviously, if $\bI\subseteq\bI'$, then $\pam(\fR_\bp,\bI,\bA)\le\pam(\fR_\bp,\bI',\bA)$.
 Especially, the number $\pam(\fR_\bp,\bA)=\pam(\fR_\bp,\bJ,\bA)$ is the biggest one.

 Proposition \ref{vers}, together with the theorem on the dimensions of fibres and the Chevalley
 theorem on the image of a regular mapping (cf. \cite[Exercise II.3.19]{ha} or
 \cite[Ch.\,I,\,\S\,5,\,Theorem 6]{sh}), implies the following result.

  \begin{corol}\label{parnum}
  Let $\kF_\bp$ be an $\bI$-family of vector rank $\fR_\bp$ over a variety $X$. For each $x\in X$ set
 $X_x=\setsuch{x'\in X}{\kF_\bp(x')\simeq\kF_\bp(x)}$ and denote 
 \begin{align*}
  &X_i=\setsuch{x\in X}{\dim X_x\le i},\\
   &\pam(\kF_\bp)=\max_i\set{\dim X_i-i}.
 \end{align*} 
 Then all subsets $X_x$ and $X_i$ are \emph{constructible} (i.e. finite unions of locally closed sets) and
 $\pam(\kF_\bp)\le\pam(\fR_\bp,\bI,\bA)$.
 \end{corol} 

 Note that the bases $\mD(\fR_\bp,\bI)$ of our almost versal families are \emph{projective}, especially
 \emph{complete} varieties. In the next section we shall exploit this property.

 Tame-wild dichotomy allows to establish the derived type of an algebra knowing the behaviour of the
 numbers $\pam(\fR_\bp,\bA)$. Namely, if $\fR_k=(r_{k1},r_{k2},\dots,r_{ks})$, set
 $|\fR_\bp|=\sum_{k,i}r_{ki}$. Since it is a maximal possible number of indecomposable summands
 of complexes of vector rank $\fR_\bp$, indecomposable complexes over derived tame algebra
 form at most one-parameter families, and the parameter number grows quadratically for derived 
 wild algebras, the following corollary is evident.

 \begin{corol}\label{partw}
  An algebra $\bA$ is derived tame \iff $\pam(\fR_\bp,\bA)\le|\fR_\bp|$ for all vector ranks $\fR_\bp$.
 \end{corol} 

 An important case is that of families of \emph{free} modules, i.e. such that $\kF_k(x)\simeq a_k\bA$ for
 some integer $k$. Namely, let  $\fR(\bA)=\row as$ (we do not suppose that $\bA$ is Morita reduced). 
 For every vector $\fB=(b_m,b_{m+1},\dots,b_n)$ we set $\fB\fA=(b_m\fA,b_{m+1}\fA,\dots,b_n\fA)$
 and write $\mD(\fB,\bI)$,  $\pam(\fB,\bI,\bA)$, etc. instead of, respectively, $\mD(\fB\fA,\bI)$,
 $\pam(\fB\fA,\bI,\bA)$, etc. When we are interested in the \emph{asymptotical} behaviour
 of parameter numbers (it is enough, for instance, to establish the derived type), we can restrict our
 considerations by free complexes only. Indeed, for a vector $\fR=\row rs$, denote
 \begin{align*}
 \big[\fR/\fA\big]&=\max\setsuch{b}{ba_i\le r_i \,\text{ for all }\, i},\\
 \big]\fR/\fA\big[&=\min\setsuch{b}{ba_i\ge r_i \,\text{ for all }\, i}.
 \end{align*}
 If $\fR_\bp=(\fR_m,\fR_{m+1},\dots,\fR_n)$, set 
 \begin{align*}
  \fB&=(b_m,b_{m+1},\dots,b_n),\quad \text{where }\, b_k=\big]\fR_k/\fA\big[ ;\\
  \fB'&=(b'_m,b'_{m+1},\dots,b'_n),\quad \text{where }\, b'_k=\big]\fR_k/\fA\big[ .
 \end{align*} 
 Then, obviously, 
 $$ 
   \pam(\fB'\fA,\bI,\bA)\le \pam(\fR_\bp,\bI,\bA)\le \pam(\fB,\bI,\bA).
 $$ 
 Especially, Corollary \ref{partw} can be reformulated as follows.

  \begin{corol}\label{freetw}
  An algebra $\bA$ is derived tame \iff $\pam(\fB,\bA)\le|\fB|\dim\bA$ for every sequence
 $\fB=(b_m,b_{m+1},\dots,b_n)$.
 \end{corol}  
 
 \section{Families of algebras. Semi-continuity}
 \label{s5}

  A (flat) \emph{family of algebras} over an algebraic variety $X$ is a sheaf $\kA$ of $\kO_X$-algebras,
 which is coherent and flat (thus locally free) as a sheaf of $\kO_X$-modules. For such a family and every
 sequence $\fB=(b_m,b_{m+1},\dots,b_n)$ one can define the function $\pam(\fB,\kA,x)=\pam(\fB,\kA(x))$.
 (Recall that here $b_k$ denote the ranks of free modules in a free complex.)
 Our main result is the upper semi-continuity of these functions. 

  \begin{theorem}\label{semi}
 Let $\kA$ be a flat family of finite dimensional algebras over an algebraic variety $X$.
  For every vector $\fB=(b_m,b_{m+1},\dots,b_n)$ the function $\,\pam(\fB,\kA,x)$
 is upper semi-continuous, i.e. all sets
 $$X_j=\setsuch{x\in X}{\pam(\fB,\kA,x)\ge j}$$
 are closed.
 \end{theorem} 
 \begin{proof}
  We may assume that $X$ is irreducible. Let $\bK$ be the field of rational functions on $X$. We consider it as a
 constant sheaf on $X$. Set $\bJ=\rad(\kA\*_{\kO_X}\bK)$ and $\kJ=\bR\cap\kA$. It is a sheaf of nilpotent
 ideals. Moreover, if $\xi$ is the generic point of $X$, the factor algebra $\kA(\xi)/\kJ(\xi)$ is semisimple. Hence
 there is an open set $U\subseteq X$ such that $\kA(x)/\kJ(x)$ is semisimple, thus $\kJ(x)=\rad\kA(x)$ for every $x\in U$.
 Therefore $\pam(\fB,\kA,x)=\pam(\fB,\kJ(x),\kA(x))$ for $x\in U$; so $X_j=X_j(\kJ)\cup X'_j$, where
 $$ 
   X_j(\kJ)=\setsuch{x\in X}{\pam(\fB,\kJ(x),\kA(x))\ge j}
 $$ 
 and $X'=X\setminus U$ is a closed subset in $X$. Using noetherian induction, we may suppose that $X_j'$ is closed,
 so we only have to prove that $X_j(\kJ)$ is closed too.

  Consider the locally free sheaf $\kH=\bigoplus_{k=m+1}^n \hom(b_k\kA,b_{k-1}\kJ)$ and the projective
 space bundle $\mP(\kH)$ \cite[Section II.7]{ha}.  Every point $h\in\mP(\kH)$ defines a set of homomorphisms
 $h_k:b_k\kA(x)\to b_{k-1}\kJ(x)$ (up to a homothety), where $x$ is the image of $h$ in $X$, and the points $h$
 such that $h_kh_{k+1}=0$ form a closed subset $\mD\subseteq \mP(\kH)$. We denote by $\pi$ the restriction onto
 $\mD$ of the projection $\mP(\kH)\to X$; it is a projective, hence closed mapping. Moreover, for every point $x\in X$
 the fibre $\pi^{-1}(x)$ is isomorphic to $\mD(\fB,\kA(x),\kJ(x))$. Consider also the group variety $\kG$ over $X$:
 $\kG=\prod_{k=m}^n\mathrm{GL}_{b_k}(\kA)$. There is a natural action of $\kG$ on $\mD$ over $X$, and 
 the sets $\mD_i=\setsuch{z\in\mD}{\dim\kG z\le i}$ are closed in $\mD$. Therefore the sets $Z_i=\pi(\mD_i)$ are 
 closed in $X$, as well as $Z_{ij}=\setsuch{x\in Z_i}{\dim\pi^{-1}(x)\ge i+j}$. But $X_j(\kJ)=\bigcup_iZ_{ij}$,
 thus it is also a closed set.
 \end{proof} 

 Taking into consideration Corollary \ref{freetw}, we obtain

  \begin{corol}\label{topen}
  For a family of algebras $\kA$ over $X$ denote 
 \begin{align*}
   X_{\mathrm{tame}}&=\setsuch{x\in X}{\kA(x)\ \text{\emph is derived tame}},\\
   X_{\mathrm{wild}}&=\setsuch{x\in X}{\kA(x)\ \text{\emph is derived wild}}  .
 \end{align*} 
 Then $X_{\mathrm{tame}}$ is a countable intersection of open subsets and
 $X_{\mathrm{wild}}$ is a countable union of closed subsets.
 \end{corol} 

 The following conjecture seems very plausible, though even its analogue for usual tame algebras has not yet been
 proved.

  \begin{conj}
  For any (flat) family of algebras over an algebraic variety $X$ the set $X_{\mathrm{tame}}$ is open.
 \end{conj} 

 Recall that an algebra $\bA$ is said to be a (flat) \emph{degeneration} of an algebra $\bB$, and $\bB$ is said to
 be a (flat) \emph{deformation} of $\bA$, if there is a (flat) family of algebras $\kA$ over an algebraic variety $X$
 and a point $p\in X$ such that $\kA(x)\simeq\bB$ for all $x\ne p$, while $\kA(p)\simeq\bA$. One easily verifies
 that we can always assume $X$ to be a \emph{non-singular curve}. Corollary \ref{topen} obviously implies
  \begin{corol}\label{def}
  Suppose that an algebra $\bA$ is a flat degeneration of an algebra $\bB$. If $\,\bB$ is derived wild,
 so is $\bA$. If $\,\bA$ is derived tame, so is $\bB$.
 \end{corol} 

 If we consider non-flat families, the situation can completely change. The reason is that the dimension is 
 no more constant in these families. That is why it can happen that such a ``degeneration'' of a derived wild algebra
 may become derived tame, as the following example due to Br\"ustle \cite{br} shows.

  \begin{exam}\label{bru}
 There is a (non-flat) family of algebras $\kA$ over an affine line $\mA^1$ such that all of them except $\kA(0)$
 are isomorphic to the derived wild algebra $\bB$ given by the quiver with relations
 $$ 
    \xymatrix{
   && \bullet \\
   \bp \ar[r]^{\al} &\bp \ar[r]^{\be_1} &\bp \ar[u]_{\ga_1} \ar[d]^{\ga_2}
   & \bp \ar[l]_{\be_2} & {\be_1\al=0,} \\
   && \bullet
  }
 $$ 
 while $\kA(0)$ is isomorphic to the derived tame algebra $\bA$ given by the quiver with relations
 \begin{equation}\label{quiv}
 \vcenter{   \xymatrix{
  & & \bp \\
  \bp \ar[r]^{\al} &\bp \ar[ur]^{\xi_1} \ar[r]^{\be_1} & \bp\ar[u]_{\ga_1} \ar[d]_{\ga_2}
  & \bp \ar[l]_{\be_2} \ar[dl]^{\xi_2} & {\be_1\al=\ga_1\be_1=\ga_2\be_2=0.}\\
 & & \bp
  } }
 \end{equation}
 Namely, one has to define $\kA(\la)$ as the factor algebra of the path algebra of the quiver as in
 \eqref{quiv}, but with the relations $\be_1\al=0,\ \ga_1\be_1=\la\xi_1,\ \ga_2\be_2=\la\xi_2$. 
 Note that $\dim\bA=16$ and $\dim\bB=15$, which shows that this family is not flat.
 \end{exam} 

 Actually, in such a situation the following result always holds.

  \begin{prop}\label{fact}
  Let $\kA$ be a family (not necessarily flat) of algebras over a non-singular curve $X$ such that
 $\kA(x)\simeq\bB$ for all $x\ne p$, where $p$ is a fixed point, while $\kA(p)\simeq\bA$. Then
 there is a flat family $\kB$ over $X$ such that $\kB(x)\simeq\bB$ for all
 $x\ne p$ and $\kB(p)\simeq\bA/I$ for some ideal $I$.
 \end{prop} 
 \begin{proof}
 Note that the restriction of $\kA$ onto $U=X\setminus \set p$ is flat, since $\dim\kA(x)$ is constant there.
 Let $n=\dim\bB$, $\Ga$ be the quiver of the algebra $\bB$ and $\bG=\fK\Ga$ be the path algebra 
 of $\Ga$. Consider the Grassmannian $\Gr(n,\bG)$, i.e. the variety of  subspaces of codimension $n$
 of $\bG$. The ideals form a closed subset $\alg=\alg(n,\bG)\subset\Gr(n,\bG)$. The restriction of the
 canonical vector bundle $\kV$ over the Grassmannian onto $\alg$ is a sheaf of ideals in
 $\kG=\bG\*\kO_{\alg}$, and the factor $\kF=\kG/\kV$ is
 a universal family of factor algebras of $\bG$ of dimension $n$. Therefore there is a morphism
 $\phi:U\to\alg$ such that the restriction of $\kA$ onto $U$ is isomorphic to
 $\phi^*(\kF)$. Since $\alg$ is projective and $X$ is non-singular, $\phi$ can be continued
 to a morphism $\psi:X\to\alg$. Let $\kB=\psi^*(\kF)$;
 it is a flat family of algebras over $X$. Moreover, $\kB$ coincides with $\kA$ outside $p$. Since both
 of them are coherent sheaves on a non-singular curve and $\kB$ is locally free, it means that
 $\kB\simeq\kA/\kT$, where $\kT$ is the torsion part of $\kA$, and  $\kB(p)\simeq\kA(p)/\kT(p)$.
 \end{proof} 
 
  \begin{corol}
  If a degeneration of a derived wild algebra is derived tame, the latter has a derived wild factor algebra.
 \end{corol} 

 In the Br\"ustle's example \ref{bru}, to obtain a derived wild factor algebra of $\bA$, one has to add the
 relation $\xi_1\al=0$, which obviously holds in $\bB$.

 By the way, as a factor algebra of a tame algebra is obviously tame (which is no more true for
 derived tame algebras!), we get the following corollary (cf. also \cite{cb,dg}). 

  \begin{corol}
  Any deformation (not necessarily flat) of a tame algebra is tame. Any degeneration of a wild algebra
 is wild.
 \end{corol}

\end{document}